\input amstex
\documentstyle{amsppt}

\NoRunningHeads \magnification=1200 \TagsOnRight \NoBlackBoxes
\hsize=5.3in \vsize=7.2in \hoffset= 0in \baselineskip=12pt
\topmatter
\title  Closed symmetric 2-differentials of the 1st kind
\endtitle
\rightheadtext{Closed symmetric 2-differentials of the 1st kind}
\author Fedor Bogomolov*
 Bruno De Oliveira**
\endauthor
\thanks
* Partially supported by NSF grant DMS-1001662 and by AG Laboratory GU-HSE grant RF government ag. 1111.G34.31.0023.
** Partially supported by the NSF grant DMS-0707097. The second author thanks the financial support  by the Courant Institute of Mathematical Sciences during 2011/12.
\endthanks
\affil
* Courant Institute for Mathematical Sciences and Laboratory of Algebraic Geometry, GU-HSE **University of Miami
\\\\
To the memory of Andrey Todorov, friend and colleague
\endaffil
\address
Fedor Bogomolov Courant Institute for Mathematical Sciences,
New York University\\
Bruno De Oliveira University of Miami
\endaddress
\email bogomolo{\@}CIMS.NYU.EDU bdeolive@math.miami.edu
\endemail

\endtopmatter
\document

\noindent {\bf Abstract:} A closed symmetric differential of the 1st kind is a differential that locally is the product of closed holomorphic 1-forms. We  show that closed symmetric 2-differentials of the 1st kind on a projective manifold $X$ come from maps of $X$ to cyclic or dihedral quotients of Abelian varieties  and that their presence implies that the fundamental group of $X$ (case of rank 2) or of the complement $X\setminus E$ of a divisor $E$ with negative  properties (case of rank 1) contains  subgroup of finite index with infinite abelianization. Other results include: i) the identification of the differential operator characterizing closed symmetric 2-differentials on surfaces (which provides in this case a connection to flat Riemannian metrics) and ii) projective manifolds $X$ having symmetric 2-differentials $w$ that are the product of two closed meromorphic 1-forms are irregular, in fact if $w$ is not of the 1st kind (which can happen), then $X$ has a fibration $f:X \to C$ over a curve of genus $\ge 1$.

\head {0. Introduction}\endhead

\

\

Our purpose is continue the study of the properties of symmetric tensors on algebraic manifolds, with a special emphasis on the case of algebraic surfaces.
This work deals with closed symmetric differentials which are a natural generalization of closed holomorphic 1-forms to symmetric differentials of higher degrees. The presence of closed holomorphic 1-forms imposes topological restrictions on the manifold, our ultimate goal is to establish a similar but weaker connection between the existence closed symmetric differentials and the topology of the manifold. 

\

A closed symmetric differential is a symmetric differential which can be decomposed as a product of closed holomorphic 1-forms on a neighborhood of some point of the manifold. If a differential has local decompositions of this type around every point of the manifold, then the  closed differential is said to be of the 1st kind. A feature of closed symmetric differentials is that they are connected to webs, possibly degenerate, on the complex manifold not matter what the dimension is while non closed symmetric differentials do not have necessarily this connection if the manifold has dimension greater than 2. Closed symmetric differentials are not necessarily of the 1st kind for degrees greater than 1. In section 2.1 we provide  examples for each of the 3 causes of this failure. One of the causes of this failure  has a striking  manifestation, some  closed symmetric differentials (holomorphic) only allow decompositions into products of closed 1-forms  if some of the 1-forms  have essential singularities.

\

In this article we study the degree 2 case  which is the first interesting case to consider and has some special features such as: 1) the local decompositions of a  differential as a product of closed holomorphic 1-forms have  rigidity properties for all degrees but for degree 2 this rigidity has the maximum strength, i.e. the 1-forms in the decomposition are unique up to a multiplicative constant; 2) symmetric 2-differentials on complex manifolds  are analogous to  Riemannian metrics in differential geometry.  Using feature 2) we identify in theorem 2.1 the differential operator that characterizes closed symmetric 2-differentials on surfaces. This differential operator is just the natural translation of the Gaussian curvature operator to our case. We observe that a closed symmetric 2-differential on a surface is the direct analogue of the notion of a flat Riemannian metric on a real surface.

\

The  rigidity of the local decompositions of a closed symmetric 2-differential as product of closed holomorphic 1-forms gives for symmetric 2-differentials of 1st kind $w$  on a projective manifold $X$ a dual pair $(\Bbb C_w,\Bbb C^*_w)$ of  local systems of rank 1 on $X$ or on an unramified double cover of $X$ (depending on whether $w$ is split or non split). Moreover, the 2-differential $w$ of the 1st kind  can be decomposed as $w=\phi_1\phi_2$ the product of  twisted  closed 1-forms $\phi_i$ ( $\phi_1$ twisted by $\Bbb C_w$ and $\phi_2$ twisted by $\Bbb C_w^*$). The work Beauville, Green-Lazarsfeld and Simpson on the cohomology loci (see references) plays a key role in establishing that the local system $\Bbb C_w$ must be torsion.  Once established that the local system $\Bbb C_w$ is torsion, one can work further to obtain the  theorem describing the geometric origins and features of non-degenerate  symmetric 2-differentials $w$ of the 1st kind, where $w$ being non-degenerate means that $w$ defines at some point $x\in X$ 2 distinct directions in $T_xX$:

\proclaim {Theorem 3.2} Let $X$ be a smooth projective manifold with $ w\in H^0(X,S^2\Omega^1_X)$ a nontrivial non-degenerate closed differential of the $1^{st}$ kind. Then:

i) $X$ has a holomorphic map to a cyclic or dihedral quotient of an abelian variety from which the symmetric differential $w$ is induced from. More precisely,  there is a commutative diagram

$$\CD
X' @> a_{X'}>> Alb(X') \\  @VfVV   @VqVV
\\ X @>a>> Alb(X')/G \\ \endCD$$

 \noindent and $f^*w=a_{X'}^*\omega$ with $\omega \in H^0(Alb(X'),S^2\Omega^1_{Alb(X')})^{G}$ where $f:X' \to X$ is an unramified $G$-Galois covering and $a_{X'}:X'\to Alb(X')$ the Albanese map. The group $G$ is  $\Bbb Z_m$ if $w$ is split and $D_{2m}$ if $w$ is non split.
\smallskip

ii) $\pi_1(X)$ is infinite, more precisely $\exists \Gamma \vartriangleleft \pi_1(X)$ such that $\pi_1(X)/\Gamma$ is finite cyclic or dihedral and its abelianization, $\Gamma/[\Gamma,\Gamma]$, is an infinite group.
\endproclaim

\

In section 2.3 we prove that if a closed 2-differential is of the 1st kind outside of codimension 2, then it is of the 1st kind everywhere. We also show that the locus where  a closed 2-differential fails to be of the 1st kind is contained in the divisorial part of the degeneracy locus, i.e. the locus of all points where the 2-differential fails to define two distinct hyperplanes on the tangent space.  

\

In our last result, theorem 3.3, we  describe  the geometry associated to the class of closed symmetric 2-differentials that are the product of two closed meromorphic 1-forms. There some points of interest in this result: one is  that this class contains symmetric differentials that are not of the 1st kind; another is that even though we are allowing the closed 1-forms $\phi_i$ in the decomposition $w=\phi_1\phi_2$ to be not holomorphic, one still obtains non-triviality of the regularity. A priori  allowing the $\phi_i$ to be not holomorphic could invalidate any non triviality results on the fundamental group (closed meromorphic 1-differentials can exist in simply connected manifolds), but in fact we show that the decomposition $w=\phi_1\phi_2$ implies the existence of a fibration over a curve of genus $\ge 1$ which implies a large fundamental group.

\proclaim {Theorem 3.3} Let $X$ be a smooth projective manifold and $w\in H^0(X,S^2\Omega^1_X)$ be a
closed differential of rank 2 with a decomposition:

$$w=\phi_1\phi_2  \text { }\text { }\text { }\text { } \text {with $\phi_i \in H^0(X,\Omega^1_{X,cl}(*))$} \tag 3.7$$

\noindent where $\Omega^1_{X,cl}(*)$ is the sheaf of closed meromorphic 1-forms. Then the Albanese dimension of $X$ $\ge 2$ and either:

\smallskip

1) $w=\phi_1\phi_2$ with $\phi_i\in H^0(X,\Omega^1_{X})$, or
\smallskip
2) $X$ has a map to a curve of genus $\ge 1$, $f:X\to C$ and $w=(f^*\varphi+u)f^*\mu$, with $f^*\varphi+u$ non-holomorphic and where $u\in H^0(X,\Omega^1_X)$, $\varphi\in H^0(C,\Omega^1_C(*))$  $\mu\in H^0(C,\Omega^1_C)$.

\endproclaim

\newpage

\head {1. Preliminaries}\endhead

\

 A symmetric differential of degree $m$, $w\in H^0(X,S^m\Omega^1_X)$, defines at each point $x \in X$  an homogeneous polynomial of degree $m$ on the tangent
space $T_xX$. If $X$ is a surface, then $w$ defines at each tangent space, $T_xX$,
$m$, not necessarily distinct, lines through the origin. Around a general point on $X$, 
one obtains $k\le m$ integrable distributions of lines giving a k-web,  i.e. a collection of k foliations.
On higher dimensions this is no longer necessarily
the case, since the pointwise splitting of $w(x)$ into linear factors  might not hold and even
if such splitting occurs the distributions of hyperplanes in $T_xX$ defined by $w$
might not be integrable.

\

Let $X$ be a complex manifold of dimension $n$. The $\Bbb P^{n-1}$-bundle $\Bbb P(\Omega_X^1)$ over $X$, $\pi:\Bbb P(\Omega_X^1) \to X$, and its tautological line bundles $\Cal O_{\Bbb P(\Omega_X^1)}(m)$ are intimately connected to the theory of symmetric differentials. There is in particular a natural bijection between $H^0(X,S^m\Omega^1_X)$ and $H^0(\Bbb P(\Omega_X^1),\Cal O_{\Bbb P(\Omega_X^1)}(m))$. To a symmetric differential $w\in H^0(X,S^m\Omega^1_X)$ on $X$ one can associate an hypersurface:

$$Z_w \subset \Bbb P(\Omega_X^1) \tag 1.1$$

\noindent such that $Z_w\cap \pi^{-1}(x)$ is an hypersurface of degree m ($Z$ can also be viewed as  the zero locus of the section of $\Cal O_{\Bbb P(\Omega_X^1)}(m)$ corresponding to $w$).

\smallskip

The hypersurface $Z_w$ can be reducible and non reduced. The irreducible components of $Z_w$ are called horizontal if they dominate $X$ via the map $\pi$ and  vertical otherwise. Hence:

$$Z_w=Z_{w,h}+Z_{w,v} \tag 1.2$$

\noindent with $Z_{w,h}$ and $Z_{w,v}$ the union of  respectively the horizontal and the vertical irreducible components.

\

\proclaim {Definition 1.1} A symmetric differential $w \in H^0(X,S^m\Omega^1_X)$
on a smooth  complex manifold $X$ is said to be:

a) split if either one of the following equivalent
statements holds:
\smallskip

\noindent i) $w=\phi_1...\phi_m$ with $\phi_i$ meromorphic 1-differentials.

\smallskip

\noindent  ii) $Z_{w,h}$ is the union of $m$ irreducible components.

\

b) split at $x$ if there is a neighborhood of $x$ on which $w$ splits (or equivalently $w(x)\in S^m\Omega^1_{X,x}$ is a product of linear forms).
\endproclaim

\noindent  Gauss lemma implies that if $w$ is split at $x$, then there is a neighborhood $U_x$ of $x$ where $w$ is the product of holomorphic 1-forms. A split symmetric differential is therefore locally the product of holomorphic 1-forms but the converse does not necessarily hold
(e.g. $Z_w$ is an unramified cover of $X$, with degree $>1$). If $w$ splits, then 
 $w=\mu_1\otimes...\otimes\mu_m$ with $\mu_i\in H^0(X,\Omega_X^1\otimes L_i)$ where the $L_i$  line bundles on $X$ with $\prod L_i=\Cal O(-D)$ with $\pi^*D=Z_{w,v}$.
 
 \

 The following fact will be used later, any  $w\in H^0(X,S^2\Omega^1_X)$ that splits at the general point (always holds if $X$ is a surface) has associated to it a canonical generically 2-1 covering (possibly ramified)  $s_w:X' \to X$ for which $s_w^*w$ is split and $X'$ is smooth.

\proclaim {Definition 1.2} On a smooth complex manifold $X$ a symmetric differential $w \in H^0(X,S^m\Omega^1_X)$
 that  splits at the general point is said to:
\smallskip
i) have rank $k$, $rank(w)=k$, if at a general point there are k distinct hyperplanes in $P(\Omega^1_{X,x})$ defined by $w(x)$, i.e.
there are $k$ distinct foliations defined by $w$ near the general point (the foliations are defined by the 1-forms $\phi_i$ in a local decomposition $w|_{U_x}=\phi_1...\phi_m$).

\smallskip

ii) be degenerate at $x$ if  $w(x)=0$ or if the number of distinct hyperplanes in $P(\Omega^1_{X,x})$ defined by $w(x)$ is less than the rank of $w$. The locus consisting of the union of all points where $w$ is degenerate is called the degeneracy locus  of $w$, $D_w$.
\endproclaim

 The degeneracy locus $D_w$ of a symmetric differential $w\in H^0(X,S^m\Omega^1_X)$ on a surface is the discriminant divisor of  $w$ which is defined locally where $w|_U=a_m(dz_1)^m+a_{m-1}(dz_1)^{m-1}dz_2+...+a_0(dz_2)^m$ by the discriminant of $w|_U$ seen as a polynomial in $\Cal O(U)[dz_1,dz_2]$. As a set $D_w$ is the the union of the points $x\in X$ such that $Z_w\cap \pi^{-1}(x)$ has multiple points.

\

\

\head {2. Closed  2-differentials and differentials of the 1st kind}\endhead

\

\

\subhead {2.1 General concepts}\endsubhead

\

\

\proclaim {Definition 2.1}  A symmetric differential $w\in H^0(X,S^m\Omega^1_X)$ on a smooth complex manifold $X$ is
 said to have:
 \smallskip
 
 1) a holomorphic (meromorphic) exact decomposition if:

$$\text { }\text { }\text { }\text { }\text { }\text { }\text { }\text { }\text { }\text { }\text { }\text { }\text { }\text { }\text { }\text { }\text { }\text { }\text { }\text { }w=df_1...df_m
\text { }\text { }\text { }\text { }\text { }\text { }\text { }\text { }\text { }\text { } f_i\in \Cal O(X) \text { }(f_i\in \Cal M(X))$$ 

\smallskip

2) a holomorphic (meromorphic) exact decomposition at $x\in X$ if there is a neighborhood $U_x$ of $x$ where $w|_{U_x}$ 
has a holomorphic (meromorphic) exact decomposition.

\smallskip

3) a split closed  decomposition if  $w|_U=\phi_1...\phi_m$ with $\phi_i$ closed holomorphic 1-differentials  on a Zariski open $U$ (the $\phi_i$ are not necessarily meromorphic on $X$).\endproclaim

\proclaim {Definition 2.2}  A symmetric differential $w \in H^0(X,S^m\Omega^1_X)$
   on a smooth complex manifold $X$ is said to be:
   \smallskip
   i) closed if $w$ has a holomorphic exact decomposition at a general point of $X$.
   \smallskip
   ii) closed of the $1^{st}$ kind if $w$ has a holomorphic exact decomposition at all points of $X$.

   \endproclaim

\

\noindent If $w$ is a symmetric differential  on a surface $X$, then  there are always holomorphic functions $f_i$ and $f$ on a neighborhood of any
general point of $X$, such that  $w=fdf_1...df_m$ holds. The condition of $w$ being closed
asks for the existence of functions $f_i$ such that $f$ can be made constant. For degree $2$ this condition
can be seen as a flatness curvature type condition on $w$ (see the next section).

\

For degree $m=1$ the classes of closed and closed of the 1st kind differentials coincide. A symmetric differential of degree 1, i.e. a holomorphic 1-form, which is closed in the sense of definition 2.2 is also closed in the usual sense due to the principle of analytic continuation.  Poincare's lemma  
implies that a closed 1-form  must be locally exact, i.e. of the $1^{st}$ kind in the sense of definition 2.2. For degrees $m\ge 2$ the two classes no longer coincide. There are 3 consecutive levels of possible the failure of a closed symmetric differential $w$ to be of  the 1st kind at $x\in X$, which will be illustrated by examples below.

\

The first level of failure  of a closed symmetric differential $w$ to be of  the 1st kind at $x\in X$ is the failing of $w$ to  split at $x$ (1st kind must be  split at every point by definition).

\

\noindent Example: (non-split) Let $z_1$ be a holomorphic  coordinate of $\Bbb C^n$  and $f\in \Cal O(\Bbb C^n)$,  set $w=z_1(dz_1)^2-(df)^2$ . The differential is non split at all points in $\{z_1=0\}$ but it is closed since any point $y \in X\setminus \{z_1=0\}$  has a  neighborhood $U_y$ where $\sqrt {z_1}$ exists and hence $w$ has a holomorphic exact decomposition $w|_{U_y}=d(\frac {2}{3}z_1^{\frac {3}{2}}+f)d(\frac {2}{3}z_1^{\frac {3}{2}}-f)$.

\

If the differential is locally split at $x$, then  the 2nd layer of failure is due to  monodromy in the factors in the exact decomposition (not of the foliations) around the locus where it fails to be of the 1st kind.

\

\noindent Example (monodromy of the exact decomposition): Let $B\subset \Bbb C^2$ be  a sufficiently small open ball about the origin where $1+z_2$ is invertible. Consider $w=(1+z_2)^\alpha dz_1d[z_1(1+z_2)]$.  The differential $w$ has an exact decomposition at a point if and only if we can decompose $(1+z_2)^\alpha$ as a product of functions of $z_1$ and $z_1(1+z_2)$ (see section 2.2).  Any such decomposition is unique up multiplicative constants (see 2.3). At points away from $\{z_1=0\}$ we have the decomposition $(1+z_2)^\alpha=z_1^{-\alpha}[z_1(1+z_2)]^\alpha$, but the functions involved are multivalued in neighborhoods of points in $\{z_1=0\}$. In fact this monodromy is infinite if $\alpha\not \in \Bbb Q$, meaning that even after finite ramified coverings the symmetric differential would not have an exact decomposition along a divisor.

\smallskip

If the differential is both locally split at $x$ and no monodromy occurs, then $w$ has  a split closed decomposition at $x$ and the 3rd level of failure is due to the singularities of the the 1-differentials on the decomposition.

\

\noindent Example: (singularities) This example shows that even essential singularities can occur, $w= e^\frac {z_2}{1+z1z_2}dz_1d[z_1(1+z_1z_2)]$. The 1-differentials in the split closed decomposition are unique up to constants, as it will shown in section 2.3, and the constants will cancel each other so in fact the decomposition is unique and has the form 

$$w=e^\frac {z_2}{1+z1z_2} dz_1d[z_1(1+z_1z_2)]=e^{-\frac {1}{z_1}}dz_1e^{\frac {1}{z_1(1+z_1z_2)}}d[z_1(1+z_1z_2)]$$

\noindent with essential singularities occurring on the closed 1-forms at $\{z_1=0\}$.

\

\noindent Next are some cases and examples of globally defined closed symmetric differentials.

\

\noindent Example:  The basic examples of  global closed
symmetric differentials  on compact Kahler manifolds are holomorphic 1-forms and their products
$w=\mu_1...\mu_m$ with $\mu_i\in H^0(X,\Omega^1_X)$. The condition of compactness is  essential t
o obtain the closed property. In the compact surface case
one does not need the Kahler condition, the holomorphic 1-forms are  closed by a direct
application of the Stoke's theorem.

\

On abelian varieties every symmetric differential is a linear combination over $\Bbb C$ of
closed symmetric differentials. 

\

On the case  of curves the space of symmetric differentials of degree
$m$ is equal to $H^0(C,(\Omega^1_C)^{\otimes m})=H^0(C,mK_C)$ and they are all closed. One of its main themes of this article 
is the decomposition of a closed symmetric differential as a product 
of  closed 1-differentials with torsion coefficients, i.e. sections of $\Omega^1_{C}\otimes \Cal O_{\chi_i}$ with 
$\Cal O_{\chi_i}$ a flat line bundle. We proceed to consider the special case of symmetric 
differentials on curves. Any symmetric differential $w\in H^0(C,(\Omega^1_C)^{\otimes m})$  is defined modulo
an invertible constant by the zero divisor of  $w$, $(w)_0$. Hence if there is a splitting of $(w)_0$ into 
collections of $(2g-2)$-points, then the splitting provides modulo constants an unique  decomposition of $w$ as a 
product of twisted 1-differentials, $w=\prod_{i=1}^m \phi_i$, where $\phi_i\in H^0(C,\Omega^1_C\otimes \Cal O_{\chi_i})$  and  
$\bigotimes_{i=1}^m \Cal O_{\chi_i}= \Cal O$. Thus we have a finite number of representations of $w$ as a product of 
1-differentials with torsion coefficients, namely there are exactly $((2g-2)m)!/((2g-2)!)^m $ such representations (modulo constants)
 if we add multiplicities. However, if we want to represent $w$ as a product of untwisted 1-differentials, i.e holomorphic 1-forms, then 
 such representation does not exists for a generic symmetric differential $w$.

\

\noindent Example (Bo-De11): All symmetric differentials $w$ of rank 1
on a projective manifold are  closed.

\

\noindent Example: let $g:Y \to X$ be an unramified covering of $X$ of
degree $m$. Then the norm for the map $g$ of an holomorphic 1-form
$\mu \in H^0(Y,\Omega^1_Y)$, $n_g(\mu)\in H^0(X,S^m\Omega^1_X)$ is locally exact.
Locally $n_g(\mu)$ is defined $n_g(\mu)|_U=df_1...df_m$ where
$df_i=\mu|_{U_i}$, $g^{-1}(U)=\prod_{i=1}^m U_i$ ($U$ sufficiently small and
the $U_i$ are biholomorphic to $U$, so the $f_i$ can also be view on $U$).

\

\noindent Example:  There are examples of global closed  symmetric differentials that are not of the fist kind (i.e. there are some points where the symmetric differential fails to have a local holomorphic exact decomposition at). These examples, that have rank 2, can be construct via (holomorphic) products of 2 closed meromorphic 1-differentials, $w=\phi_1\phi_2$ with $\phi_i\in H^0(X,\Omega^1_{cl}(*))$ and $\phi_1\wedge \phi_2$ is somewhere non-vanishing. In section 3.3, we show that if we have such a product, then $X$ has a fibration $f:X \to C$ onto a smooth 
curve of genus $\ge 1$ and $w=(f^*\varphi+u)f^*\mu$, where $u\in H^0(X,\Omega^1_X)$,
$\varphi\in H^0(C,\Omega^1_C(*))$ and $\mu\in H^0(C,\Omega^1_C)$ ($\Omega^1_C(*)$ is the sheaf of meromorphic 1-differentials on $C$). 

\

\

\subhead {2.2  Differential operator for closed 2-differentials on surfaces}\endsubhead

\

\

Poincare lemma states that locally exact symmetric differentials $w$ of degree 1
 are the solutions of the first order differential equation $dw=0$,
where $d$ is the exterior derivative. For higher degrees, we saw in the last section, the 
illustrative examples of how a symmetric differential that has almost everywhere a local holomorphic exact 
decomposition can fail to have it  at points where the differential is degenerate, i.e. in the discriminant locus $D_w$. The reasons for this failure, e.g. monodromy about the the divisor $D_w$ (giving a local cohomological obstruction to the existence of a holomorphic exact decomposition), can not be detected via a differential operator. A differential operator 
can only be expected to detect the existence of local holomorphic decompositions where the symmetric differential is non degenerate, which by the way is enough to guarantee that the symmetric differential is closed. The result of this section states that  the property of a symmetric 2-differential being  closed is indeed determined by a differential operator. We expect the same to happen for higher degrees (see the end of the section).

\

A symmetric 2-differential $w\in H^0(X,S^2\Omega^1_X)$  is locally given by:

$$w(z)|_U=a_{11}(z)dz_1^2+a_{12}(z)dz_1dz_2+a_{22}(z)dz_2^2 \tag 2.1$$

\noindent with $a_{ij}(z)\in \Cal O(U)$ and can formally be considered as a  (degenerate) "complex metric" on $X$. This perspective illustrates
once more the distinction between rank 1 and rank 2. Only the case of rank 2 benefits from this perspective, since rank 1 would 
correspond to an everywhere degenerate metric. The reasoning that follows connecting the property of being closed to flatness  concerns rank 2 alone. The case of rank 1 is distinct and related to the case of 1-forms,  we have that locally near a general point $w|_{U}=f(z)(dz_1)^2$ and $w$ is closed if and only if $df\wedge dz_1=0$. Moreover, global arguments give that all symmetric differentials of rank 1 defined on a compact Kahler manifold are closed, see [BodeO11].

\

The Gaussian curvature $R$ operator on 2-dimensional real manifolds acts on
sections of $S^2(T^{\Bbb R}_xM)^*$ representing metrics and sends them to functions. We
call $R^{\Bbb C}$ the natural "complexification" of this operator (i.e. replace the
$\frac {\partial}{\partial x_i}$ by $\frac {\partial}{\partial z_i}$) which gives a map:

$$R^{\Bbb C}:H^0(X,S^2\Omega^1_X)\to \Cal M(X)$$

\noindent Recall that associated with the symmetric differential of rank 2 $w$ one has $\det (w)$ the  section of $2K_X$, given locally  by $\det (w)(z)=a_{11}(z)a_{22}(z)-\frac {1}{4}a_{12}(z)^2$.

\proclaim {Theorem 2.1} Let $X$ be a smooth complex surface and  
 $w\in H^0(X,S^2\Omega^1_X)$ a symmetric differential of rank 2. Set $P_2w=\det(w)^2R^{\Bbb C}$. Then  the nonlinear differential operator 
 
 $$P_2: H^0(X,S^2\Omega^1_X) \to H^0(X,4K_X) \tag 2.2$$
 
 \noindent is such that  $w$ is closed
 if and only if $P_2w=0$. Moreover, $P_2w=0$ implies is of the 1st kind on $X\setminus D_w$.
 \endproclaim

\demo {Proof} The complexified Gaussian curvature $R^{\Bbb C}$ operator applied to the symmetric 2-differential $w$ which is locally expressed in the form (2.1)  gives according to the Brioschi formula:
$$R^{\Bbb C}w|_U=\frac {1}{det(w|_U)^2} [\vmatrix  -\frac {1}{2}(\frac {\partial^2a_{11}}{\partial z_2^2}+
\frac {\partial^2a_{22}}{\partial z_1^2}-\frac {\partial^2a_{12}}{\partial z_1\partial z_2})&\frac {1}{2}\frac {\partial a_{11}}{\partial z_1}&
\frac {1}{2}(\frac {\partial a_{12}}{\partial z_1}-
\frac {\partial a_{11}}{\partial z_2})
\\
\frac {1}{2}(\frac {\partial a_{12}}{\partial z_2}-
\frac {\partial a_{22}}{\partial z_1})&a_{11}&\frac {1}{2}a_{12}
\\
\frac {1}{2}\frac {\partial a_{22}}{\partial z_2}&\frac {1}{2}a_{12}&a_{22}\endvmatrix
- \vmatrix 0&\frac {1}{2}\frac {\partial a_{11}}{\partial z_2}&\frac {1}{2}
\frac {\partial a_{22}}{\partial z_1}
\\
\frac {1}{2}\frac {\partial a_{11}}{\partial z_2}&a_{11}&\frac{1}{2}a_{12}
\\
\frac {1}{2}\frac {\partial a_{22}}{\partial z_1}&\frac{1}{2}a_{12}&a_{22}
\endvmatrix]$$

\noindent Globally one obtains a meromorphic function $R^{\Bbb C}w$ whose poles come from the zeros of $det(w)^2$ and (2.2) follows.

\

Every point outside of the discriminant locus of $w$, $x\in X\setminus D_w$,  has an open neighborhood  $U_x$ with $w|_{U_x}=\mu_1\mu_2$ with $\mu_i$ nowhere vanishing holomorphic 1-forms  on $U_x$.  The existence of local non vanishing holomorphic integrating factors for  nowhere vanishing holomorphic 1-forms implies that  after possibly shrinking once more $U_x$ one has:

$$w|_{U_x}=f(u)du_1du_2 \tag 2.3$$

 \noindent with $f\in \Cal O(U_x)$ and  $(u_1,u_2)$ a holomorphic coordinate chart of $U_x$. 
 
 \
 
 The "Gaussian curvature"  for "complex metric"  in the form (2.3) is given by:

$$R^{\Bbb C}w|_{U_x}=-\frac {2}{f}\frac{\partial ^2{log f}}{\partial {u_1}\partial {u_2}} \tag 2.4$$

\noindent ($f$ is non vanishing). 

\

The condition $R^{\Bbb C}w|_{U_x}=0$ (i.e. $\frac{\partial ^2{log f}}{\partial {u_1}\partial {u_2}}=0$) is equivalent to   $f(u)=f_1(u_1)f(u_2)$
on some ball $B_x$ centered at $x$. Hence $w$ has a holomorphic  exact decomposition on $B_x$ 

$$w(u)=dF_1(u_1)dF_2(u_2)$$

\noindent where the $F_i$ are the primitives of the $f_i$.  Hence  the condition that $P_2w=0$ is equivalent to the existence of a locally exact decomposition
of $w$ at every point in $X\setminus D_w$, which implies the theorem.

\enddemo

\

\proclaim {Question} Is there a differential operator $P_m$, $m>2$,  generalizing $P_2$ and characterizing closed
 symmetric  m-differential on surface $X$?
\endproclaim

The following is a sketch of an approach  to show that such operators or better said a family of such differential 
operators do exist. Let $w$ be a symmetric m-differential on a surface $X$. Locally on an sufficiently 
small open neighborhood $U_x$ of a generic point $x\in X$ the symmetric tensor $w$ is given by 
the product 

$$w|_{U_x}=\mu_1^{m_1}...\mu_k^{m_k}$$ 

\noindent with $\mu_i\in H^0(U_x,\Omega^1_X)$ and  $\mu_i\wedge \mu_j$ nowhere vanishing 
for $i\neq j$. In this case $w$ defines a nonsingular k-web $\Cal W_w$, i.e. a family of  k distinct foliations $\Cal F_i$ on $X$ which are 
pairwise transversal and smooth on $U_x$. Let $\{z_i\}_{i=1,...,k}$ 
with $z_i\in \Cal O(U_x)$ be a set of local functions such that: 1) $dz_i$ are nowhere vanishing and 
2) whose level sets are the leaves of the foliations $\Cal F_i$ on $U_x$. For such a collection $\{z_i\}_{i=1,...,k}$
one gets the decomposition $w|_{U_x}= f \prod_{i=1}^k dz_i^{m_i}$ with $f\in \Cal O(U_x)$. 

\

Let us consider the germ a nonsingular m-web $\Cal W$ at a point $x\in X$  defined by the symmetric 
differential $dz_1...dz_m$ with  $\{z_i\}_{i=1,...,m}$ a collection as above. Any germ $w_x$ of a symmetric 
differential $w$ at $x$ with $\Cal W_{w_x}=\Cal W$ can be written in the form 

$$w_x=f\prod_{i=1}^mdz_i \tag 2.5$$

\noindent with $f\in \Cal O_x$. The jet of n-th order of a symmetric m-differential defining $\Cal W$
is determined by $(n+2)(n+1)/2$ coefficients (from the Taylor series of $f$), giving the dimension of $J^n(\Cal W,x)$, the space 
of n-th order jets of symmetric m-differentials defining the  m-web germ $\Cal W$. We denote by $J^n(\Cal W,x)_{cl}\subset J^n(\Cal W,x)$
the space of n-th order jets of closed m-differentials defining the m-web $\Cal W$. 
The symmetric differential $w_x$ is closed if and only if the function $f$  is of the form $f=\prod_{i=1}^m f_i(z_i)$. The n-th order jets of each function $f_i$ (determined by $n+1$ coefficients) are involved in the n-th order jet of the product function $f$ and the constant term of the product $f$ imposes only one condition of the n-th order jets of the $f_i$. Consequently, $J^n(\Cal W,x)_{cl}$ is a subvariety of $J^n(\Cal W,x)$ of dimension at most $mn+1$ and hence it is a proper subvariety once $n>2m-3$.

\

Let $J^n(m,x)$ ($J^n(m,x)_{cl}$) be the set of n-th order jets of symmetric (closed) m-differentials 
at $x$. From the previous discussion it follows that once $n>2m-3$ the closure $\bar J^n(m,x)_{cl}$ 
of set $J^n(m,x)_{cl}$ is a proper affine subvariety of $J^n(m,x)$. On $J^n(m,x)$ we have the natural 
action of $J^nAut(\Bbb B^2)$ which is the group of $n$-th order jets of  holomorphic automorphisms 
of the ball the 2-ball $(\Bbb B^2)$. The group $J^n Aut(\Bbb B^2)$ is an algebraic group which is  
finite-dimensional nilpotent extension of $GL(2)$. The subvariety $\bar J^n(m,x)_{cl}$ is naturally 
invariant under the above action of $J^n Aut(\Bbb B^2)$.

\

Let $F$ be a regular function on $J^n(m,x)$ which vanishes on $\bar J^n(m,x)_{cl}$
and satisfies $F(gw_x)= \chi^N(g) F(w_x)$  $\forall w_x\in J^n(m,x)$ and $\forall g\in J^n Aut(\Bbb B^2)$
with $\chi :J^n Aut(\Bbb B^2) \to GL(2)\to C^*$ natural projection and $N\in \Bbb N^+$ 
(the function $F$ with the latter property will be called a semi-invariant function relative 
to the group action). Standard invariant theory gives that there is a finite set of semi-invariant
regular functions $F$ generating the ideal of semi-invariant functions vanishing on $\bar J^n(m,x)_{cl}$.
Thus the function $F$ defines a nonlinear map of vector bundles over $X$: 

$$F' : J^n(m,x)\to NK_X$$ 

\noindent with $F'$ mapping the (nonlinear) sub-fibration $\bar J^n(m,x)_{cl}$
into  zero section of $NK_X$. 
The map $F'$ on the n-th order jets of symmetric m-differential induces a differential operator of order $n$

$$D_F: H^0(X,S^m\Omega^1_X) \to H^0(X,NK_X)$$

\noindent which is trivial on the closed symmetric m-differentials
which define nonsingular m-webs. Our previous discussion stating that
 if $n>2m-3$, then  $\bar J^n(m,x)_{cl}\subsetneq J^n(m,x)$ implies that the n-th order differential operators 
 just described will vanish on the closed symmetric m-differentials but will be nontrivial
 on generic symmetric m-differentials. Note that for the case of $m=2$ this approach gives
 that we need to go to jets of order 2 to obtain a differential operator which vanishes on closed 
 but not on generic symmetric differentials (this matches result in theorem 2.1). The  construction 
 just described raises up many interesting questions. The most fundamental questions are clearly
 how to find such semi-invariant functions $F$ and operators $D_F$ naturally
for arbitary $m$ and what are the properties of such operators.

\

\

\subhead {2.3 1st kind, local systems and global decompositions}\endsubhead

\

In this section we describe the analytical and topological objects (respectively twisted holomorphic closed 1-forms and local systems) that can be associated to a closed symmetric 2-differential of the 1st kind and establish two basic facts about the locus where  a closed symmetric 2-differential fails to be of the 1st kind. In the next section we will determine the properties and geometric consequences of these objects. We remind the reader that the case of interest concerns 2-differentials of rank 2, the case of rank 1 will be mentioned just in passing (for full details see [BoDeO11]).

\

The local holomorphic exact decompositions of a closed symmetric differential $w$ 
$$w_{|_U}=(df_1)^{m_1}...(df_k)^{m_k} \tag 2.6$$

\noindent with $df_i\wedge df_j\not \equiv 0$, have rigidity properties. The strength of the rigidity is
dependent on a notion coming from the theory of webs, the  abelian rank of the k-web associated to the differential $w$ on $U$ (for the notion of abelian rank see for example [ChGr78]).
In the case of interest, i.e. differentials of degree 2,  then the decomposition (2.6) has the strongest form of rigidity, i.e. the $f_i$ are unique up to a multiplicative constant and an additive constant  (in the case of rank 2 this is a manifestation of the triviality of the abelian rank of the 2-web and will be directly explained below). 

\

\proclaim {Remark 2.2} i) A rank 2 symmetric 2-differential of the 1st kind is always locally of the form (2.6) (by definition).

\smallskip

ii)  The case of a rank 1 symmetric 2-differential of the 1st kind is distinct, there might be points where $w$ can not be locally written in the form (2.6). 
The example to have in mind is $w=z_1(dz_1)^2$, which can not be put in the form (2.6) on any open set intersecting $\{z_1=0\}$ (due to the absence of a well defined square root). In general,  let $(w)_0=\sum_il_iD_i$, be the irreducible decomposition of $(w)_0$. The local decomposition of type (2.6) exist in a neighborhood of every point outside of the divisor:

$$E_w=\sum_{j\in \{i| m\nmid l_i\}}D_j$$
 
The presence of $E_w$ has an impact on the topological properties that can be derived from the presence of a rank 1 symmetric differential (see section 3.1 and [BoDeO11] for full details). 
\endproclaim
\

\proclaim {Lemma 2.3} Let $X$ be a smooth complex manifold and $w$ a symmetric 2-differential having  a  closed  decomposition, i.e. $w|_U=\phi_1\phi_2$, with $\phi_i\in H^0(U,\Omega^1_{_{\text {cl}}})$ for some Zariski open set $U\subset X$. 
\smallskip
i) If $w$ is of rank 2, then the closed decomposition of $w$ is unique, up to  multiplication by constants.  More precisely, if $V\subset X$ is any connected open subset of $U$ and $w|_V=\psi_1\psi_2$ is another closed decompositions of $w|_V$, then up to a reordering of the $\psi_i$ we have:

$$\text { }\text { }\text { }\text { }\text { }\text { }\text { }\text { }\text { }\text { }\text { }\text { }\psi_i=c_i\phi_i|_V \text { }\text { }\text { }\text { }\text { }\text { }\text { }\text { }\text { }\text { }\text { } c_2=c_1^{-1}\in \Bbb C^*$$

ii) If $w$ is of rank 1, then the closed  decompositions of $w$ are not unique. However,  any two decompositions of the form $w=(\phi)^{ 2}=(\psi)^{ 2}$ with $\phi,\psi \in H^0(V,\Omega^1_{_{\text {cl}}})$, $V\subset X$  open, are also unique up to  multiplication by $\pm 1$.
\endproclaim

\demo {Proof} The  case of rank 1 is clear. Consider the case where $w$ is of rank 2. Let $x\in X$  and $U_{x}$ be an open neighborhood of $x$ where $\phi_i|_{U_x}=df_i$ and $\psi_i|_{U_x}=dg_i$ with $f_i,g_i\in \Cal O(U_x)$. The differential $w$ being of rank 2 implies that we could have chosen $x\in X$  such that $df_1(x)\wedge df_2(x)\neq 0$, i.e. $f_1$ and $f_2$ can be viewed as local holomorphic coordinates
around $x$. \par

After reordering, one can make $df_i\wedge dg_i\equiv 0$ for $i=1,2$.
The relation $df_i\wedge dg_i\equiv 0$ implies that near $x$ the $g_i$ is a function of $f_i$, $g_i=g_i(f_i)$. Hence
$dg_1dg_2=df_1df_2$ implies that $g'_1(f_1)g'_2(f_2)=1$ must hold, which can only happen if $g'_1(f_1)$ and $g'_2(f_2)$ are
nontrivial constant functions (since $df_1(x)\wedge df_2(x)\neq 0$). So $\psi_i=c_i\phi_i$ on a neighborhood of $x$ and hence on the 
whole $V$ via the principle of analytic continuation.
\enddemo

\

\proclaim {Remark/Notation 2.4} A symmetric 2-differential of the 1st kind and rank 2 on a complex manifold $X$ is either split or there is  an associated  unramified double cover of $X$

$$s_w:X' \to X \tag 2.7$$

\noindent such that $s_w^*w$ splits. The differential not being split is equivalent to the impossibility to make the local orderings of the two foliations associated to $w$ consistent on the whole $X$. Given a open covering $\Cal U$ of $X$ with local orderings, one gets from the transition of the orderings in the intersections a 1-cocycle in $H^1(\Cal U,S_2)$ and its associated representation  $\rho:\pi_1(X) \to S_2$.  The covering $s_w$ is the regular cover associated with the representation $\rho$.\endproclaim

\

The next proposition shows that associated to a split 2-differential $w$ of the 1st kind one has a dual pair of local systems and a decomposition of $w$ as a product of twisted closed holomorphic 1-forms. This decomposition will be used in section 3 to derive the global geometric/topological properties of $X$. If the  differential $w$ is non split, then the remark above tell us that the pair $(X',s_w^*w)$  has this decomposition from which we derive again the geometric/topological properties of $X$.

\proclaim {Proposition 2.5} Let $X$ be a smooth complex manifold and $w\in H^0(X,S^2\Omega^1_X)$ be split and closed of the 1st kind. 
\smallskip
i) (rank 1) $w$ has a decomposition on $X\setminus E_w$, $w|_{E_w}=\phi^2$ with $\phi \in H^0(X\setminus E_w,\Omega^1_{_{\text {cl}}}\otimes_{\Bbb C}
\Bbb C_{w})$, where the divisor $E_w$ is as in remark 2.2, $\Bbb C_w$ is a local system of rank 1 associated to a 1-cocycle with values in $\Bbb Z_2$.

\smallskip
ii) (rank 2) $w$ has a decomposition:

$$w=\phi_1\phi_2 \tag 2.8$$

\noindent with $(\phi_1,\phi_2)\in H^0(X,\Omega^1_{_{\text {cl}}}\otimes_{\Bbb C}
(\Bbb C_{w}\oplus \Bbb C^*_{w}))$ where $\Bbb C_w$ is a local system of rank 1 whose isomorphism class is uniquely 
 determined (up to its dual) by $w$.
\endproclaim

\demo {Proof}  Case of rank 1, see [BoDeO11]. Case of rank 2. Let $\Cal U=\{U_i\}_{i\in I}$  be a 
covering of $X$ by holomorphic balls such that their intersection are contractible. This covering can be chosen such that $w|_{U_i}=df_{1i}df_{2i}$ with $f_{1i},f_{2i}\in \Cal O(U_i)$ (1st kind) and $df_{1i}\wedge df_{1j}\equiv 0$ on all nonempty intersections $U_i\cap U_j$ (split). Applying Lemma 2.3, it follows that:

 $$\text { }\text { }\text { }\text { }\text { }\text { }\text { }\text { }\text { }\text { }\text { }\text { }\text { }\text { }\text { }df_{k i}=c_{k,ij}df_{k j}  \text { }\text { }\text { }\text { }\text { }\text { }\text { }\text { }  \text{ and } c_{2,ij}=c_{1,ij}^{-1}\tag 2.9$$

\noindent where $k=1,2$ and $c_{k,ij}\in \Bbb C^*$. The 2 collections $\{c_{k,ij}\}$ for $k=1,2$ are 1-cocycles in $Z^1(\Cal U,\Bbb C^*)$. It follows from (2.9) that:

$$(\phi_1,\phi_2):=(\{df_{1 i}\}_{i\in I}, \{df_{2 i}\}_{i\in I})\in H^0(X,\Omega^1_{_{\text {cl}}}\otimes (\Bbb C_{w}\oplus\Bbb C^*_{w} ))$$

\noindent  where $\Bbb C_{w}$ is the local system of rank 1 determined by the 1-cocycle $\{c_{1,ij}\}$.\smallskip

The pair of local systems $(\Bbb C_{w},\Bbb C^*_{w} )$ is uniquely determined up to isomorphism (and order, of course) by the 2-differential of the 1st kind $w$.  Let $\Cal V=\{V_k\}$ be another  cover of $X$ for which the $w|_{V_k}$ are holomorphically exact decomposable with $w|_{V_k}=dg_{1k}dg_{2k}$ and let $\Cal W=\{W_r\}$ a Leray cover relative to locally constant sheaves of $X$ such that each $W_r$ is such that $W_r\subset U_{i(r)}$ and $W_r\subset V_{k(r)}$. The induced decompositions $w|_{W_r}=df_{1i(r)}|_{W_r}df_{2i(r)}|_{W_r}$ and $w|_{W_r}=dg_{1k(r)}|_{W_r}dg_{2k(r)}|_{W_r}$ following the argument above can be used to derive two pairs of 1-cocycles with values in $\Bbb C^*$ relative to the cover $\Cal W$. Since by  lemma 2.3 the following $df_{\alpha i(r)}|_{W_r}=c_{\alpha,r}dg_{\alpha k(r)}|_{W_r}$ holds for $\alpha=1,2$, it follows that the 1-cocycles are cohomologous and hence the pairs of associated local systems are isomorphic.
\enddemo

\

\noindent Example:  The non-triviality of the local systems associated to a split closed symmetric 2-differential does occur. Let $C$  be a curve with an involution without fixed points. The surface $Y=(C\times C)/\Bbb Z_2$ where $\Bbb Z_2$
acts diagonally ($Y=C/\Bbb Z_2 \times C/\Bbb Z_2$) is such that  $q(C\times C)-q(Y)= g(C)-1$. So if $g(C)\ge 3$ there are 2 anti-invariant
holomorphic 1-forms $\mu_1$ and $\mu_2$ on $C\times C$, one coming from each factor. The product $\mu_1\mu_2$ is $\Bbb Z_2$-invariant and hence induces a 2-differential $w$ on $Y$ which is closed, rank 2 and split (the 2 foliations of $w$ correspond to the 2 natural fibrations on $C/\Bbb Z_2 \times C/\Bbb Z_2$). The symmetric 2-differential $w$ is the product of the twisted closed 1-forms on $Y$ coming from the $\mu_i$ (and can not be written as the product of two untwisted 1-forms).

\

The next lemma shows  that the locus where a closed differential fails to be of the 1st kind has no isolated points and hence is of pure codimension 1.
\

 \proclaim {Lemma 2.6} Let $X$ be a complex manifold and $w\in H^0(X,S^2\Omega^1_X)$ be  closed of the 1st kind outside of codimension 2. Then $w$ is closed of the 1st kind on $X$.
\endproclaim

\demo {Proof}  Let $Z\subset X$ be a  locus of  codimension at least 2 containing all the points of $X$ where $w$ fails to have holomorphic exact decomposition at. Pick any $x\in Z$ and $B_x\subset X$ a ball centered at $x$. Let $\Cal U=\{U_i\}$ be an open covering
of $B_x\setminus Z$ on which:

$$w|_{U_i}=df_{1i}df_{2i} \tag 2.10$$

\noindent with $f_{1i}, f_{2i}\in \Cal O(U_i)$. Since $\pi_1(B_x\setminus Z)=0$
we can order for each $i$ the functions $f_{ki}$ such that on the non-empty intersections  $U_{ij}$ $df_{ki}\wedge df_{kj}=0$.
By the lemma 2.3 on the intersections $U_{ij}$ one has:

$$df_{ki}=c_{k,ij}df_{kj}$$

\noindent  where $\{c_{k,ij}\} \in Z^1(\Cal U,\Bbb C^*)$ (for k=1 and 2).  These cocycles must be coboundaries since
$B_x\setminus Z$ is simply connected. Hence for each $k$ there is a 0-cochain $\{c_{k,i}\}$
with values in $\Bbb C^*$ such that $d(c_{k,i}f_{k,i})=d(c_{k,j}f_{k,j})$ on $U_{ij}$.  The collections
$\{d(c_{k,i}f_{ki})\}$ glue to give two closed  1-forms $\mu_1, \mu_2 \in H^0(B_x\setminus Z,\Omega^1_{_{\text {cl}}})$ such that:

 $$w|_{B_x\setminus Z}=\mu_1\mu_2$$
 
 \noindent Again since
$\pi_1(B_x\setminus Z)=0$ it follows that the forms $\mu_i$ are actually exact, i.e.
$\mu_i=df_i$ with $f_1, f_2 \in \Cal O(B_x\setminus Z)$.  Hartog's extension theorem
gives the holomorphic extensions  $\bar f_1, \bar f_2 \in \Cal O(B_x)$ of respectively $f_1$ and $f_2$
from which   $w$ gets the holomorphic exact decomposition at $x$, $w|_{B_x}=d\bar f_1d\bar f_2$. Hence there are no pints where $w$ fails to be of the first kind.
\enddemo

\

\proclaim {Proposition 2.7} Let $X$ be a smooth complex manifold and  
 $w\in H^0(X,S^2\Omega^1_X)$ a closed  differential of rank 2. Then the locus where $w$ fails to be of the 1st kind is contained in the  degeneracy locus $D_w$.
\endproclaim

\demo {Proof} Since $X\setminus D_w$ is connected and the locus where $w$ is of the 1st kind is open and nonempty, we just need to show that the set where $w$ is of the 1st kind is also closed in $X\setminus D_w$.\par

We want to prove that if $x\in X\setminus D_w$ is such that  all open balls $B_x$  centered at $x$  have a point $y \in B_x$ where $w$ has a holomorphic exact decomposition at, then $w$ also has a holomorphic exact decomposition at $x$. Since $x\in X\setminus D_w$, then $w$ is split at $x$ and one has a ball $B_x$ centered at $x$ where 

$$w|_{B_x}=\mu_1\mu_2$$ 

\noindent with $\mu_i\in H^0(B_x,\Omega^1_X)$. Using the hypothesis of the claim there is a point $y\in B_x$  at which $w$ has a holomorphic  exact decomposition at, 

$$w|_{B_y}=dh_1dh_2 \tag 2.11$$

\noindent with $h_i\in \Cal O(B_y)$, where $B_y\subset B_x$ is an open ball centered at $y$. 
It follows from (2.11) that, after reordering if necessary, one has $\mu_i|_{B_y}=g_idh_i$ with 
$g_i\in \Cal M(B_y)$, hence $\mu_i\wedge d\mu_i=0$ on $B_y$. This in turn implies the $\mu_i$ are completely integrable on the whole $B_x$. Since $\mu_i(x)\neq 0$ and the $\mu_i$ are completely integrable  one can shrink $B_x$ so that the $\mu_i$ have first integrals $u_i\in \Cal O(B_x)$, i.e. $\mu_i=f_idu_i$ with $f_i,u_i\in \Cal O(B_x)$.  Since $\mu_1\wedge \mu_2$ is nowhere vanishing on $B_x$, again by shrinking $B_x$ we can assume $u_1$ and $u_2$ are two coordinates of a holomorphic chart $(u_1,..,u_n)$ for $B_x$. So one has

$$w|_{B_x}=h(u)du_1du_2$$

\noindent (with $h(u)=f_1f_2 \in \Cal O(B_x)$). In the ball $B_x$ there is  a point $y$ where $w$ has a holomorphic exact decomposition at, still describe this decomposition by (2.1). From $dh_1dh_2= hdu_1du_2$  on $B_y$ it follows that $h_i=h_i(u_i)$ and hence 

$$h(u)|_{B_y}=h_1'(u_1)h'_2(u_2)$$

\noindent This in turn implies that the function $h$ depends only on $u_1$ and $u_2$, $h=h(u_1,u_2)$, on the whole $B_x$ and $h(u_1,u_2)=h_1'(u_1)h'_2(u_2)$ on $B_y$. The latter condition, as discussed in the proof of theorem 2.1, implies that that $h$ satisfies the differential equation $\frac{\partial ^2{log h}}{\partial {u_1}\partial {u_2}}=0$ on $B_y$. By the identity principle  

$$\frac{\partial ^2{\log h(u_1,u_2)}}{\partial {u_1}\partial {u_2}}=0 \tag 2.12$$

\noindent  holds on $B_x$. The equation 2.12 as seen in the proof of theorem 2.1 implies that  $h=h_1(u_1)h_2(u_2)$ on a neighborhood of $x$ and $w$ therefore has a holomorphic exact decomposition at $x$.

\enddemo

\

\

\head {3. Global geometric properties}\endhead

\

\

In the first part of this section we  determine the topological restrictions and geometric features that are implied by closed symmetric 2-differentials of the 1st kind. Here we are mainly interested in the properties of the fundamental group and the existence of varieties (and maps into them) from which the differentials would be induced. The case of  2-differentials of   rank 1 follows from a previous work by the authors [BoDeO11], hence the focus lies in the rank 2 case. In the last subsection, we describe the geometry of a very natural of class of closed 2-differentials (not necessarily of the 1st kind), the class consisting of products of two closed meromorphic 1-forms.

\

\

\subhead {3.1 1st kind of rank 1}\endsubhead

\

\

In [BoDeO11] the authors studied symmetric differentials of rank 1 of all degrees on projective manifolds. We present here, for the sake of completeness, the statement of the result concerning the case of interest, i.e. degree 2, with a small modification concerning the condition of being of 1st kind  plus a few remarks. The main result of that paper gives for the case of interest:

\proclaim {Theorem 3.1} [BoDeO11] Let $X$ be a smooth projective manifold and $w\in H^0(X,S^2\Omega^1_X)$
a nontrivial differential of rank 1.  Then:

\smallskip i) $w$ is closed on $X$ and closed of the 1st kind outside of codimension 2. 

\smallskip

ii)There is a cover of $X$ which is generically 2 to 1 $g:X' \to X$ such that $g^*w=\mu^{\otimes 2}$ with $\mu\in H^0(X',\Omega^1_{X'})$.

 \smallskip 
 
 iii)  There is a holomorphic map  from $X$ to a $\Bbb Z_2$-quotient  of an abelian variety with isolated singularities
 $a_w:X \to A_w/\Bbb Z_2$, such that $w={a_w}^*(u)$ and $u\in H^0(A_w/\Bbb Z_2,
 S_{\text {orb}}^m\Omega^1_{A_w/\Bbb Z_2})$.

 \smallskip iv) There is a 2-negative divisor $E\subset (w)_0  \subset X$ such that $\pi_1(X\setminus E)$ is infinite. More precisely, $\pi_1(X\setminus E)$
 has a normal subgroup $\Gamma$ of finite index for which $\pi_1(X\setminus E)/\Gamma$ is cyclic  and its abelianization, $\Gamma/[\Gamma,\Gamma]$, is an infinite group.
\endproclaim

\

\noindent $\bold {Remarks}$: 1) A divisor $D$ is said to be 2-negative if for all smooth surfaces $S\subset X$, the divisor $D\cap S$  of $S$ is negative.

\smallskip

2) item iii) states that geometrically a symmetric 2-differential of rank 1 (of 1st kind or not) comes from an orbifold symmetric differential on a $\Bbb Z_2$-quotient of an abelian variety.  If $w=\mu^{\otimes 2}$ with $\mu\in H^0(X,\Omega^1_{X})$, then  $a_X$ is the natural $a_w$, where $a_X$ is the Albanese morphism.

\smallskip

3) The divisor $E$ lies inside the ramification divisor of the covering map $g$ in ii) (which in turn lies inside the union of the irreducible components of $(w)_0$ with odd multiplicity). 

\smallskip

4) For the case of closed 2-differentials of the 1st kind of rank 2 we will see below that one obtains topological conditions for the whole complex manifold $X$, while in this case (of rank 1) the conditions are for the complement $X\setminus E$. The reason for this distinction is the fact that for 2-differentials of the first kind of rank 1 the rigidity of the holomorphic exact decompositions at the points in $E$ is weakened (outside of $E$ one can write $w=(df)^2$ and such decomposition is unique at to multiplication by $\pm 1$ but for example $w=z^3(dz)^2$ has no natural exact decomposition that is unique). 

\smallskip 5) There are projective manifolds with  a closed symmetric 2-differential of the 1st kind and rank 1 which are  simply connected (see [BoDeO11]).

\

\

\subhead {3.2 1st kind of rank 2}\endsubhead

\

\

In order to extract the geometrical/topological properties associated to closed symmetric 2-differentials of the 1st kind we are going to take full advantage of their global decomposition of  as products of twisted closed 1-differentials. As seen in proposition 2.5, if $w$ is split $w=\phi_1\phi_2$, with $(\phi_1,\phi_2)\in H^0(X,\Omega^1_{_{\text {cl}}}\otimes_{\Bbb C} (\Bbb C_{w}\oplus \Bbb C^*_{w}))$ where $(\Bbb C_{w}\oplus \Bbb C^*_{w})$ is the dual pair of local systems associated to $w$ via the closed decomposition. \par

There are 2 elements in the above decomposition with a geometric meaning. One is the dual pair of local systems $(\Bbb C_w,\Bbb C_w^*)$ giving us a pair of dual characters of the fundamental group. The other are the twisted closed 1-differentials $\phi_i$ which topologically define elements of 1st cohomology group of $X$ with coefficients in the local system $\Bbb C_w$ or $\Bbb C_w^*$ and geometrically define special  type of foliations, e.g. if $\Bbb C_w\otimes \Cal O$ is non-torsion the foliations are algebraic (it follows from the work Beauville, Green-Lazarsfeld and Simpson on the cohomology loci (see references)) and if torsion the foliations share the same properties as foliations defined by global holomorphic 1-forms. A key ingredient of next result is that if $X$ is projective, then  the isomorphism class of the local system $\Bbb C_w$ must be torsion. Hence in the projective case, if $w$ is split , then up to finite cyclic unramified covers, rank 2 differentials of the $1^{st}$ kind are just products of closed 1-differentials.
\

\proclaim {Theorem 3.2} Let $X$ be a smooth projective manifold with $ w\in H^0(X,S^2\Omega^1_X)$ a nontrivial rank 2 closed differential of the $1^{st}$ kind. Then:

i) $X$ has a holomorphic map to a cyclic or dihedral quotient of an abelian variety from which the symmetric differential $w$ is induced from. More precisely,  there is a commutative diagram

$$\CD
X' @> a_{X'}>> Alb(X') \\  @VfVV   @VqVV
\\ X @>a>> Alb(X')/G \\ \endCD$$

 \noindent and $f^*w=a_{X'}^*\omega$ with $\omega \in H^0(Alb(X'),S^2\Omega^1_{Alb(X')})^{G}$ where $f:X' \to X$ is an unramified $G$-Galois covering and $a_{X'}:X'\to Alb(X')$ the Albanese map. The group $G$ is  $\Bbb Z_m$ if $w$ is split and $D_{2m}$ if $w$ is non split.
\smallskip

ii) $\pi_1(X)$ is infinite, more precisely $\exists \Gamma \vartriangleleft \pi_1(X)$ such that $\pi_1(X)/\Gamma$ is finite cyclic or dihedral and its abelianization, $\Gamma/[\Gamma,\Gamma]$, is an infinite group.
\endproclaim

\demo {Proof}  We start by considering the case when $w$ is split. A split closed 2-differential of the first kind and rank 2 has as stated in proposition 2.5 the decomposition:

$$w=\phi_1\phi_2 \tag 3.1$$

\noindent with $(\phi_1,\phi_2)\in H^0(X,\Omega^1_{_{\text {cl}}}\otimes_{\Bbb C}
(\Bbb C_{w}\oplus \Bbb C^*_{w}))$ where $\Bbb C_w$ is a local system of rank 1 whose isomorphism class is uniquely determined up to its dual by $w$.  We will first show how the theorem follows if the isomorphism class of local system $\Bbb C_w$ is torsion and then prove  that the class of $\Bbb C_w$ is indeed torsion.

\

Assume: the  isomorphism class of $\Bbb C_w$ is torsion.

\

Associated to the local system $\Bbb C_w$ (whose isomorphism class is torsion) we have a finite character $\rho_w: \pi_1(x) \to \Bbb C^*$, with image a cyclic group $\Bbb Z_m$.  The unramified cyclic Galois  cover

$$f:X' \to X$$  

\noindent associated to the character $\rho_w$ is such that $f^*\Bbb C_w$ is isomorphic to the trivial local system $\Bbb C$ on $X'$. \par

First we show how to use the pullback of the decomposition (3.1)  to $X'$ to obtain the   decomposition 

$$f^*w=\mu_1\mu_2$$

\noindent with $\mu_i\in H^0(X',\Omega^1_{X'})$ (despite the $f^*\phi_i$ not  being the 1-forms $\mu_i$). Let $\Cal U=\{U_i\}_{i\in I}$ be a Leray  open covering of $X'$ relative to locally constant sheaves. On the contractible open sets $U_i$ one has $f^*\phi_{k}|_{U_i}=dg_{k,i}$ with  $g_{k,i}\in \Cal O(U_i)$ for $k=1,2$. As in the proposition 2.5, on the intersections $U_i\cap U_j$, $dg_{k,i}=b_{k,ij}dg_{k,j}$, with $b_{1,ij}=b_{2,ij}^{-1}$ and $\{b_{1,ij}\}$ the 1-cocycle  relative to $\Cal U$ giving the local system $f^*\Bbb C_w$. Since the local system $f^*\Bbb C_w$ is isomorphic to the trivial local system $\Bbb C$ on $X'$ there are 0-cochains relative to $\Cal U$ with values in $\Bbb C^*$, $\{b_{k,i}\}_{i\in I}$ ($b_{1,i}=b_{2,i}^{-1}$),   whose coboundaries are the 1-cocycles $\{b_{k,ij}\}$ for $k=1,2$. To obtain the 1-forms $\mu_k$ giving $f^*w=\mu_1\mu_2$  one untwists the collections $\{dg_{k,i}\}_{i\in I}$ using the 0-cochains $\{b_{k,i}\}_{i\in I}$ 

$$\mu_k=\{b_{k,i}dg_{k,i}\}_{i\in I}$$

\par

The $\Bbb Z_m$ action on the Galois covering space $X'$ of $X$   induces due to the universal properties of the Albanese variety of $X'$ an action on $Alb(X')$ and the Albanese morphism $a_{X'}:X' \to Alb(X')$ descends to the morphism $a:X\to Alb(X')/\Bbb Z_m$  asked in part i) of the theorem. It follows from $f^*w=\mu_1\mu_2$ with $\mu_k\in H^0(X',\Omega^1_{X'})$ that $f^*w=a_{X'}^*\omega$ with $\omega \in H^0(A,S^2\Omega^1_X)^{\Bbb Z_m}$. Finally the topological consequence, part ii), is a direct consequence of i).

\

Claim: the isomorphism  class of $\Bbb C_w$ is torsion.

\

Case: isomorphism class of $\Bbb C_w$ is non torsion but $L_w=\Bbb C_w\otimes \Cal O$ is torsion.

\

Since $X$ is a compact kahler manifold there is  an unique isomorphism class of unitary local systems giving any fixed flat line bundle. Let  $\Bbb C_u$  be a unitary local system such that $L_w\simeq\Bbb C_u\otimes \Cal O$.  From $L_w$ being torsion plus the uniqueness of the isomorphism class of unitary local systems giving the trivial line bundle it follows that the isomorphism class of $\Bbb C_u$  is also torsion. Therefore, as above, there is a finite unramified covering  $f:X' \to X$ such that  $f^*\Bbb C_u\simeq \Bbb C$ and hence  $f^*L_w\simeq \Cal O$. Note that $f^*\Bbb C_w$ is not isomorphic to the trivial local system.\par

\

Consider the pullback, $f^*w=f^*\phi_1f^*\phi_2$, of the decomposition (3.1) to $X'$. Let $\Cal U=\{U_i\}_{i\in I}$ be a Leray  open covering of $X'$ relative to locally constant sheaves where  $\{f^*\phi_{k}|_{U_i}\}_{i\in I}=\{dg_{k,i}\}_{i\in I}$ with  $g_{k,i}\in \Cal O(U_i)$ and $k=1,2$. Hence

$$f^*w|_{U_i}=dg_{1,i}dg_{2,i} \tag 3.2$$

\noindent On the intersections $U_i\cap U_j$, $dg_{k,i}=b_{k,ij}dg_{k,j}$, with  $\{b_{1,ij}\}$ the 1-cocycle  relative to $\Cal U$ giving the local system $f^*\Bbb C_w$ and $b_{1,ij}=b_{2,ij}^{-1}$. 

\

The first paragraph of this case tells us that  while the cohomology class  $[\{b_{1,ij}\}]\in H^1(X',\Cal O^*)$ is trivial, the cohomology class $[\{b_{1,ij}\}]\in H^1(X',\Bbb C^*)$ is nontrivial. Let $\{h_i\}_{i\in I}$ be the 0-cochain with values in $\Cal O^*$ whose coboundary is $\{b_{1,ij}\}$ then set

$$\mu_1=\{h_idg_{1,i}\}_{i\in I} \text { } \text { } \text { } \text { } \text {and  } \text { } \text { } \text { } \text { } \text { } \text { } \text { } \mu_2=\{h_i^{-1}dg_{2,i}\}_{i\in I}$$

\noindent By construction both collections $\{h_idg_{1,i}\}_{i\in I}$ and $\{h_i^{-1}dg_{2,i}\}_{i\in I}$ match on the intersections making $\mu_k\in H^0(X',\Omega^1_{X'})$ of $k=1,2$. Since $X'$ is compact kahler $d\mu_k=0$ and hence 

$$\{h_idg_{1,i}\}_{i\in I} =\{d\hat g_{1,i}\}_{i\in I}  \text { } \text { } \text { } \text { } \text {and  } \text { } \text { } \text { } \text { } \text { } \text { } \text { } \{h^{-1}_idg_{2,i}\}_{i\in I} =\{d\hat g_{2,i}\}_{i\in I} \tag 3.3$$

\noindent for some $\hat g_{1,i},\hat g_{2,i}\in \Cal O(U_i)$. Therefore we get for each $U_i$ 

$$f^*w|_{U_i}=d\hat g_{1,i}d\hat g_{2,i} \tag 3.4$$

\noindent and end up with two holomorphic exact decompositions of $w|_{U_i}$,  (3.2) and (3.4). Since $w$ is of rank 2, it follows from lemma 2.3  that the $h_i$ are actually constant. Hence $\{h_i\}_{i\in I}$  is 0-cochain with values in $\Bbb C^*$ but this leads to $[\{b_{1,ij}\}]\in H^1(X,\Bbb C^*)$ being trivial, a contradiction.

\

Case: $L_w$ (and  $\Bbb C_w$) is non-torsion.

\

In this case we will use the geometric properties of twisted holomorphic
differentials which were studied in [GrLa87], [Be92], [Si93] and [Ar] to understand the
cohomology locus $S_m^1(X)=\{L\in Pic^\tau(X)|\dim H^1(X,L)\ge m\}$, where $Pic^\tau(X)$ is
the variety of line bundles with trivial Chern class.\par

Consider an irreducible component $Z$ of $S^1_m(X)$ containing $L_w$  with $m=\dim H^0(X,L_w)$ and let $\alpha\in H^1(X,L_w)$ be the image of $\phi_2$ from (3.1) via the complex anti-linear isomorphism  

$$H^0(X,\Omega^1_X\otimes L_w^*) \to H^1(X,L_w) \tag 3.5$$ 

\noindent provided by conjugation. The work of Simpson [Si93] states that since the line bundle $L_w$ is non-torsion, the variety $Z$ must positive dimensional. By the construction hypothesis of $Z$  one has that  $\dim H^1(X,L)\ge \dim H^1(X,L_{w})$ (in fact equality holds due to Grauert's semi-continuity theorem) for $L$ in $Z$  hence the class $\alpha$ is preserved under small deformations of $L_w$ along $Z$. Using the work of [GrLa87] on the deformation theory of line bundles in $S^1_m(X)$, the class $\alpha$ being preserved by small deformations along $Z$ implies that 

$$\phi_2\wedge u_1=0 \tag 3.6$$

\noindent where $u_1\in H^0(X,\Omega^1_X)$ is the conjugate of  $v_1\in T_{L^*_w}Pic^\tau(X)=H^1(X,\Cal O)$ giving a 1st order deformation of
$L_w$ in $Z$ preserving $\alpha$.

\

Beauville  [Be92] obtained a Castelnuovo-De Franchis type theorem for twisted
forms from the condition $\phi_2\wedge u_1=0$. The Beauville-Castelnuovo-De Franchis
theorem states that there is a connected fibration
$f_1:X \to C_1$, $C_1$ a smooth curve, such that:

\smallskip

1) $u_1\in f_1^*H^0(C_1,\Omega^1_{C_1})$

\smallskip

2) $L_w,L_w^*\in Pic^\tau(X,f_1)$

\smallskip
\noindent where $Pic^\tau(X,f_1)$ is the subvariety of $Pic^\tau(X)$ consisting of line bundles whose
restrictions on the smooth fibers of $f_k$ are trivial. 
\noindent  The conditions (3.6) and 1) imply that the fibers $f_1$ are the leaves of the foliation defined by $\phi_2$.

\

Now repeat the previous argument replacing the line bundle $L_w$ by  the line bundle $L_w^*$ and obtain a map $f_2:X \to C_2$, 
where $C_2$ is a smooth curve, such that the fibers of $f_2$ are the leaves of the foliation defined by $\phi_1$ and $L_w,L_w^*\in Pic^\tau(X,f_2)$. 
\par
In conclusion, we have that the non-torsion line bundle $L_w$ is trivial along the smooth fibers of both fibrations $f_1:X\to C_1$ and $f_2:X \to C_2$. This can only happen if the fibrations share all fibers and hence the foliations $\phi_1$ and $\phi_2$ have the same leaves, which can not happen since $w$ has rank 2.

\

If the symmetric differential $w$ is non split, then consider the unramified double cover 

$$s_w:X'' \to X$$

\noindent described in the remark 2.4  which was built such that $(s_w)^*w$ is split. Apply the previous results for the pair $(X'',(s_w)^*w)$ to obtain an unramified $\Bbb Z_m$ covering $f:X' \to X''$ such that $(s_w\circ f')w=\mu_1\mu_2$ with the $\mu_i\in H^0(X',\Omega^1_{X'})$. The covering $f:=(s_w\circ f'):X' \to X$ is an unramified Galois cover of $X$ with Galois group $D_{2m}$ and the result follows as above.
\enddemo

\

\

\subhead {3.3 Products of closed meromorphic 1-differentials}\endsubhead

\

In this section we describe the geometry associated to closed holomorphic symmetric 2-differentials $w$ of rank 2
 which are the product of two closed meromorphic 1-forms, $w=\phi_1\phi_2$.   We will see that such $w$ might not be 
 of the 1st kind and that if that is the case, then  there must exist a fibration over a curve of genus $\ge 1$ and the Albanese dimension 
 of $X$ will (still) be greater or equal to 2. 
\

\proclaim {Theorem 3.3} Let $X$ be a smooth projective manifold and $w\in H^0(X,S^2\Omega^1_X)$ be of rank 2 with a decomposition:

$$w=\phi_1\phi_2  \text { }\text { }\text { }\text { } \text {with $\phi_i \in H^0(X,\Omega^1_{X,cl}(*))$} \tag 3.7$$

\noindent where $\Omega^1_{X,cl}(*)$ is the sheaf of closed meromorphic 1-forms. Then the Albanese dimension of $X$ is $\ge 2$ and one of the following  cases holds:

\smallskip

1) $w=\phi_1\phi_2$ with $\phi_i\in H^0(X,\Omega^1_{X})$, i.e. $w$ is of the 1st kind.
\smallskip
2) $X$ has a map to a curve $C$ of genus $\ge 1$, $f:X\to C$ and $w=(f^*\varphi+u)f^*\mu$, with $f^*\varphi+u$ non-holomorphic and where $u\in H^0(X,\Omega^1_X)$, $\varphi\in H^0(C,\Omega^1_C(*))$  $\mu\in H^0(C,\Omega^1_C)$.

\endproclaim

\demo {Proof} If the differentials $\phi_i$  have no poles, then 1) holds and the Albanese dimension is $\ge 2$ since $\phi_1\wedge \phi_2\not \equiv 0$ (rank 2). \par

From now on we assume that at least the differential $\phi_1$ has poles. First we show that we can reduce to the case where   $(\phi_i)_0\cap (\phi_i)_\infty=\emptyset$ for $i=1,2$,  where $(\phi_i)_0$ is the divisorial component of the zero locus of $\phi_i$. This follows from the next lemma and the fact that the conclusion of the theorem will hold for the pair $X$ and $w$ if it holds for the pair $X'$ and $\sigma^*w$, where $\sigma:X' \to X$ is a modification of $X$.\par

\proclaim {Lemma 3.4} Let $X$ be a smooth projective manifold and $\phi \in H^0(X,\Omega^1_{X,cl}(*))$. Then there is a modification of $X$ $\sigma:X' \to X$ such that:

$$ (\sigma^*\phi)_0 \cap (\sigma^*\phi)_\infty=\emptyset $$

\noindent and both divisors are normal
crossings divisors.
\endproclaim

\demo {Proof} Let $\sigma:  X' \to X$ be a map consisting of a finite composition of blow ups such that
the union of the divisor of poles and zeros of $\sigma^*\phi$ is a divisor with normal crossings. There is  a local chart $(z_1,...,z_n)$ near any $x\in (\sigma^*\phi)_0 \cap (\sigma^*\phi)_\infty$ such that the closed differential $\sigma^*\phi$  is of the form  

$$\sigma^*\phi|_{U_x}=\sum_{i=1}^kc_i\frac {dz_i}{z_i}+df $$

\noindent Moreover $(\sigma^*\phi)_0|_{U_x}=\bigcup_{i=l}^{dim X}\{z_i=0\}$ with $k<l$, $c_i\in \Bbb C$ and
$f\in \Cal M(U_x)$. The first observation is that the $c_i=0$ (no logarithmic pole), since $\sigma^*\phi|_{\{z_{_l}=0\}}=0$
implies the residue of $\sigma^*\phi$ about $\{z_i=0\}$, $i=1,...,k$ is 0. So $x\in (df)_0\cap (df)_\infty$ and therefore must be a point of
indeterminancy of $f$. After a finite number of further blowing ups one can resolve the indeterminancies of $f$
and the result follows (note that this lemma is not true without the assumption that $\phi$ is closed).
\enddemo

\

  Pick a connected component of the support of the polar divisor of $\phi_1$, $P\subset (\phi_1)_\infty$. Since $w=\phi_1\phi_2$ is holomorphic, then 
 
 $$P\subset Z \subset (\phi_2)_0$$ 
 
 \noindent where $Z$ is the support of the connected component of $(\phi_2)_0$ containing $P$. 
 
 \
 
The following  shows how the presence of $P$ (and then $Z$) implies the existence of a fibration with $P=Z$ as a  set theoretic fiber and $\phi_2$ is induced from the base of this fibration.\par
 
Since $Z$ is a closed analytic subvariety, $Z$ has an open neighborhood $U_Z$ such that $Z$ is a deformation retract of $U_Z$. Moreover, due to $(\phi_2)_0\cap (\phi_2)_\infty=\emptyset$, $U_Z$ can be chosen so that $\phi_2$ is holomorphic on $U_Z$. This implies that the periods of $\phi_2$ on $U_Z$ are just the periods of $\phi_2$ on $Z$ and hence they all vanish and we can integrate $\phi_2$ to get a holomorphic function. Set $\hat f(z)=\int_{z_0}^z \phi_2$ with $z_0\in Z$ and $z\in U_Z$.  By construction $Z\subset \hat f^{-1}(0)$, since $Z\subset (\phi_2)_0$. In fact, we have more:
 
 \
 
{ \bf Claim}: $Z$ is a connected component of $\hat f^{-1}(0)$.
 
 \
 
We can prove the claim by going to dimension 2. Consider a smooth surface $S=X\cap (H_1\cap ... \cap H_{\dim X-2})$, where the $H_i$ are general hyperplanes with its inclusion map $j:S\hookrightarrow X$, the induced differential $j^*w=j^*\phi_1j^*\phi_2$, $Z'=Z\cap (H_1\cap ... \cap H_{\dim X-2})$ and $P'=P\cap (H_1\cap ... \cap H_{\dim X-2})$.\par

Say $Z$ is not a connected component of  $\hat f^{-1}(0)$, then $Z'$ will not be a connected component of $\hat f^{-1}|_{S\cap U_Z}(0)$, but  $P'$ will be a connected component of $(j^*\phi_1)_\infty$ since Supp$(j^*\phi_1)_\infty=$Supp$[(\phi_1)_\infty\cap (H_1\cap ... \cap H_{\dim X-2})]$. To prove the claim, we have that  a slight of the Zariski's lemma (lemma 8.2 of [BHPV03]) gives that $Z'$ is a negative divisor, which leads to a contradiction since $P'\subset Z'$ and connected components of polar divisors of closed meromorphic 1-forms can not be negative (see below).

\
 
{Subclaim}: $P'$ is not a negative divisor.
 
 \
 
Assume $P'$ is negative. Let $P'=\sum_{r=1}^mP'_r$ be the irreducible  decomposition of $P'$ and $U'\subset S$ a sufficiently small open neighborhood of $P'$ such that $(j^*\phi_1)_\infty\cap U'=P'$.  
  \par

First we show that  the differential $j^*\phi_1$ is of the 2nd kind on $U'$, i.e. all residues along divisors on  $U'$ vanish. Recall that the residue of $j^*\phi_{1}$ along a divisor $D$ on $U'$ is $c_D=\frac {1}{2\pi i}\int_{\gamma} j^*\phi_{1}$, where $\gamma$ is a simple loop around $D$. Clearly the residues of $j^*\phi_{1}$ on $U'$ can only occur along the irreducible components $P'_r$, $r=1,...,m$, of $P'$. It follows from [Mu61] that the simple loops $\gamma_r$ around the $P'_r$, $r=1,...,m$, as elements in $H_1( {U'}\setminus P',\Bbb Z)$ satisfy the nonsingular system of m linear relations  $\sum_{j=1}^mP'_j.P'_i\gamma_i=0$, determined by the negative definite intersection matrix $[P_j.P_i]$. Hence, we have the linear system

$$\sum_{j=1}^mP'_j.P'_i\int_{\gamma_i}j^*\phi_1=0$$

\noindent which implies that all residues $c_r=\frac {1}{2\pi i}\int_{\gamma_r} j^*\phi_{1}=0$.\par

The differential $j^*\phi_1$ being of the 2nd kind on $U'$ implies that there is  a Leray open covering $\Cal U=\{U_i\}$ relative to locally constant sheaves of $U'$ such that:

$$j^*\phi_{1}|_{U_i}=dg_i \tag 3.8$$

\noindent with $g_i\in \Cal M(U_i)$ and 

 $$(j^*\phi_{1})_\infty\cap U'=\sum_{r=1}^mn_rP'_r \tag 3.9$$

\noindent with $n_r\ge2$ for all $r$. Set $P''=\sum_{r=1}^m(n_r-1)P_r$. It follows from (3.8) 
and (3.9) that $(g_i)_\infty=P''|_{U_i}$. If $\{h_i=0\}$, $h_i\in \Cal O(U_i)$,  are the equations defining the divisor $P''$, then 

$$g_i=\frac {v_i}{h_i}$$

\noindent with $v_i\in \Cal O(U_i)$  such that $v_i|_{P'_r\cap U_i} \not\equiv 0$, $r=1,...,m$. On $U_i\cap U_j$ the equality $g_i=g_j+d_{ij}$  holds with $d_{ij}\in \Bbb C$ and therefore it follows that:

$$v_i=\frac {h_i}{h_j}v_j+d_{ij}h_i \tag 3.10$$

\noindent According to (3.10) the collection $\{v_i\}$ gives a section:

$$v:=\{v_i\} \in H^0(P'',\Cal O_{P''}(P''))$$

\noindent which does not vanish identically on any irreducible component $P'_r$ of $P''$.  Hence ${P'_r}.P''\ge 0$ for all $r$ and then by linearity one has $(P'')^2\ge 0$  which contradicts the assumption that $P'$ is negative and ends the proof of the claim.

\

 We just established that $Z$ is a connected component of the level set $\hat f^{-1}(0)$. Using the compactness of $Z$ we can shrink  $U$ so  that $Z=\hat f^{-1}(0)$ as a set. The open mapping theorem and $Z$ being compact imply that there is a sufficiently small open disc $\Delta_\epsilon$ centered at $0$ such that 

$$\hat f|_{\hat f^{-1}(\Delta_\epsilon)}:\hat f^{-1}(\Delta_\epsilon)\to \Delta_\epsilon \tag 3.11$$ 

\noindent is a proper fibration onto $\Delta_\epsilon$ with $\text {Supp}\hat f^{-1}(0)=Z$. \par

The local fibration $\hat f|_{\hat f^{-1}(\Delta_\epsilon)}:\hat f^{-1}(\Delta_\epsilon)\to \Delta_\epsilon$ implies the existence of a global connected fibration $f:X \to C$ where $C$ is a smooth curve and $P=Z$ occur (set theorectically) as a fiber  and $\phi_2=f^*\mu$ with $\mu\in H^0(C,\Omega^1_C(*))$. \par 

The fact that the existence of the local fibration (3.11) implies the existence of a global fibration is a well known result. For completeness we mention a  result of this type (stronger than we need) that asserts that if $X$ has 3 connected effective divisors that are pairwise disjoint and belong to the same rational cohomology class in $H^2(X,\Bbb Q)$, then $X$ has a unique connected fibration onto a smooth curve with the divisors as fibers, see [To00].\par

The conclusion that $\phi_2=f^*\mu$ for some meromorphic 1-form $\mu\in H^0(C,\Omega^1_C(*))$ is a consequence of 

$$\phi_2\wedge f^*\eta=0 \tag 3.12$$ 

\noindent where $\eta$ is any meromorphic 1-form on $C$. The vanishing of  $\phi_2\wedge f^*\eta$ holds since by construction the fibers of $f$  are leaves of the foliation defined $\phi_2$.  It follows from (3.12) that $\phi_2=gf^*\eta$ for some $g\in \Cal M(X)$, but since both $\phi_2$ and $f^*\eta$ are closed, one has $dg\wedge f^*\eta=0$, hence $g$ is constant along the fibers of $f$ and comes from $\hat g\in \Cal M(C)$ giving $\phi_2=f^*\mu$ with $\mu=\hat g\eta$. \par

\

The fibration $f:X \to C$ has also the crucial property that the polar divisors $(\phi_1)_\infty$ and $(\phi_2)_\infty$ are contained in its  fibers.  Any connected component of $(\phi_1)_\infty$ or $(\phi_2)_\infty$, other than our previous $P$  of course, will not intersect $P=Z$, hence since $P=Z$ is a full fiber of $f$ they must be contained in a finite collection of fibers of $f$. Moreover, one has that the connected components of $(\phi_1)_\infty$ and $(\phi_2)_\infty$ occur as full fibers of $f$ as follows from the arguments that showed that $P'$ is not negative.

\

We proceed to show that $\phi_2$ is actually holomorphic and therefore $\phi_2=f^*\mu$ with $\mu\in H^0(C,\Omega^1_C)$. Suppose $(\phi_2)_\infty\neq \emptyset$, then the argument that gave the existence of the connected fibration $f:X \to C$ from $(\phi_1)_\infty\neq \emptyset$ gives that there is a connected fibration $f':X \to C'$ such that $\phi_1=f^*\mu'$ with $\mu'\in H^0(C,\Omega^1_C(*))$. According to the previous paragraph both connected fibrations $f:X\to C$ and $f:X' \to C'$   share  $P=Z$ as a fiber. This implies that the fibrations must coincide (see next paragraph). However the fibrations can not coincide since this would imply that symmetric differential $w$ would be of rank 1.\par
  Each fibration gives a holomorphic function on a neighborhood $U$ of $P$ whose level sets are the fibers. If the fibrations were distinct, then there would be fibers of one of the fibrations, contained in $U$, that would not be level sets of the holomorphic function defining the other fibration. This would give rise to non constant holomorphic functions on some fibers  which can not happen. This forces the connected fibrations  $f$ and $f'$ to coincide. 

\

What remains is to give a description of the differential $\phi_1$. The differential $\phi_1$ induces on the general fiber $F$ of $f$ a holomorphic 1-form $i^*\phi_1 \in H^0(F,\Omega^1_F)$, where $i:F \hookrightarrow X$ is the inclusion map. The differential $i^*\phi_1$ is holomorphic since all connected components of $(\phi_1)_\infty$ must be, as shown as above,  fibers of $f$. The global invariant cycle theorem by Deligne [De71] states that
if $i^*\phi_1$ remains invariant under the monodromy action of $\pi_1(C\setminus S)$, $S$ the critical values of $f$, then
there is a

$$u\in H^0(X,\Omega^1_X) \text { }\text { }\text { }\text { } \text {such that} \text { }\text { }\text { }\text { }i^*u=i^*\phi_1 \tag 3.13$$.

\noindent The invariance of $i^*\phi_1$ under the monodromy action is guaranteed by $\phi_1$ being a closed holomorphic 1-form on $f^{-1}(C\setminus S)$.
The pullback of difference $\phi_1-u$ to the general fiber of $f$ vanishes and hence  $\phi_1-u=f^*\varphi$ with $\varphi$ meromorphic differential on $C$, completing the proof of 2). The Albanese dimension in the case 2) is also $\ge 2$, the two holomorphic 1-forms $f^*\mu$ and $u$ must satisfy $f^*\mu\wedge u\neq 0$ since $w$ has rank 2.

\enddemo

\end
\Refs


\ref  \key Ar92 \by {\text { }\text { }\text { } \text
{ }\text { }\text { }\text { }\text { }\text { }D.Arapura} \paper Higgs line bundles, Green-Lazarsfeld sets, and maps of KŠhler manifolds to \text { }\text { }
\text { }\text { }\text { }\text { }\text { }\text { }\text {
}\text { }\text { }curves \jour Bull. Amer. Math. Soc. (N.S.) \vol 26 \yr 1992 \pages no. 2, 310Ð-314
\endref

\ref \key BHPV03 \by {\text { }\text { }\text { }\text { }\text {
}\text { }\text { }\text { }\text { }\text { }W.Barth,K.Hulek,C.Peters,A.Van de Ven}
\book Compact complex surfaces \publ Springer Verlag, 
Ergeb-\text { }\text { }
\text { }\text { }\text { }\text { }\text { }\text { }\text {
}\text { }nisse series vol. 4, \yr 2003
\endref

\ref  \key Be92 \by {\text { }\text { }\text { }\text { } \text
{ }\text { }\text { }\text { }\text { } A.Beauville} \paper
Annulation du $H^1$ pour les fibres en droites plats, Complex Algebraic Varieties, \text { }\text { }
\text { }\text { }\text { }\text {
}\text { }\text { }\text { }Proc. Conf., Bayreuth, 1990
 \jour Springer lecture Notes in Math \vol 1507  \yr 1992\pages 1-15
\endref

\ref  \key BoDeO06 \by {\text { }\text { }\text { }\text { } \text
{ }\text { }\text { }\text { }\text { } F.Bogomolov, B.De
Oliveira} \paper Hyperbolicity of nodal hypersurfaces \jour J.
Regne Angew. Math.  \text { }\text { }
\text { }\text { }\text { }\text {
}\text { }\text { }\text { }\vol 596 \yr 2006\pages 89-101
\endref

\ref  \key BoDeO11 \by {\text { }\text { }\text { }\text { } \text
{ }\text { }\text { }\text { }\text { } F.Bogomolov, B.De
Oliveira} \paper Symmetric differentials of rank 1 and holomorphic maps \text { }\text { }
\text { }\text { }\text { }\text {
}\text { }\text { }\text { }\text { } \text { }\jour Pure and Applied Mathematics Quarterly \vol 7 \yr 2011 \pages no. 4, 1085-1104
\endref



\ref \key ChGr78 \by  {\text { }\text { }\text { }\text { }  \text {
}\text { }\text { }\text { } \text { }S.S.Chern, P.Griffiths}
\paper Abel's theorem and webs \jour Jahresber. Deutsch. Math.-Verein.
 \vol 80 \text { }\text { }\text { }\text { }  \text {
}\text { }\text { }\text { } \text { }  \text {
}\text { } \yr {1978}  \pages  no. 1-2, 13-110
\endref



\ref  \key De71 \by {\text { }\text { }\text { }\text { }  \text {
}\text { }\text { }\text { }\text { } P. Deligne} \paper Th$\acute{e}$orie de Hodge, II
 \jour  Inst. Hautes ƒtudes Sci. Publ. Math.
\vol 40 \yr 1971 \pages 5-57
\endref



\ref \key Gr-La87 \by  {\text { }\text { }\text { }\text { }  \text {
}\text { }\text { }\text { } \text { }M.Green, R.Lazarsfeld}
\paper Deformation theory, generic vanishing theorems, and some
\text { }\text { }
\text { }\text { }\text { }\text { }\text { }\text { }\text { } \text { }\text { }
\text { }\text { }\text { }\text { }\text { }
\text { }conjectures of Enriques, Catanese and Beauville\jour
Invent. Math. \vol 90   \yr {1987}  \pages 389-407
\endref

\ref \key La04 \by {\text { }\text { }\text { }\text { }\text {
}\text { }\text { }\text { }\text { }\text { } R. Lazarsfeld}
\book Positivity in algebraic geometry \publ Springer Verlag,
Ergebnisse series vol. 48, \text { }\text { }\text { }\text {
}\text { }\text { }\text { }\text { }\text { }\text { }\text { }49\yr 2004
\endref

\ref  \key Mu61 \by {\text { }\text { }\text { }\text { }  \text {
}\text { }\text { }\text { }\text { } D.Mumford} \paper The topology of normal singularities of an algebraic surface and a criterion \text { }\text { }\text { }\text {
}\text { }\text { }\text { }\text { }\text { }\text { }\text { }for simplicity
 \jour  Inst. Hautes ƒtudes Sci. Publ. Math.
\vol 9 \yr 1961 \pages 5-22
\endref



 \ref  \key Si93 \by {\text { }\text { }\text { }\text { } \text
{ }\text { }\text { }\text { }\text { } C.Simpson} \paper
Subspaces of moduli spaces of rank one local systems
 \jour Ann. scient. Ec. Norm.  \text { }\text { }
\text { }\text { }\text { }\text {
}\text { }\text { }\text { } Sup., $4^e$ serie \vol 26 \yr 1993\pages 361-401
\endref

\ref  \key To00 \by {\text { }\text { }\text { }\text { }  \text {
}\text { }\text { }\text { }\text { } B.Totaro} \paper The topology of smooth divisors and the arithmetic of abelian varieties
 \jour Michigan \text { }\text { }
\text { }\text { }\text { }\text {
}\text { }\text { }\text { }\text { }Math. J
\vol 48 \yr 2000 \pages 611--624
\endref

 
\endRefs
